\documentclass[10pt, oneside, openright]{article}
\usepackage{epsfig}
\usepackage{graphicx, epstopdf}
\usepackage{graphpap}
\usepackage{amsmath, bm}
\usepackage{amsfonts}
\usepackage{amssymb}
\usepackage{makeidx}
\usepackage{amsmath,graphicx,url,times, amssymb, bm, epstopdf}
\usepackage{amssymb,amsmath,graphicx,times,url, mathtools}

\begin{document}

{\centering
{\Large Lyapunov Function for the Nonlinear Moog Voltage Controlled Filter}\\

\vspace{0.5in}
Technical Report AAG2021a\\
 Acoustics and Audio Group, University of Edinburgh\\

\vspace{0.5in}

Stefan Bilbao\\

\vspace{0.5in}

March 19 2021\\}

\vspace{0.5in}

In this short report, a new Lyapunov function for the Moog voltage-controlled filter is demonstrated, under zero-input conditions, and under nonlinear autonomous conditions (i.e. when parameters are not time-varying). The new definition allows for a proof of stability over the entire allowable range of parameters (cutoff frequency and resonance), and can be used as a starting point for Hamiltonian-based numerical simulation methods. \\

\newpage
\clearpage

\section{Introduction: The Moog voltage-controlled filter}
\label{introsec}

The operation of the Moog ladder filter or voltage-controlled filter (VCF) \cite{Moog65}, including nonlinear effects, is covered by H{\'e}lie \cite{Heliedafx06, Heliedafx11}. The digital emulation of the Moog VCF dates back to work by Stilson and Smith \cite{stilson96}, who considered virtual analog modeling of the Moog VCF in the linearised case, and Huovilainen, who examined the nonlinear case \cite{huovilainen_non-linear_nodate}. 

Consider the nonlinear Moog voltage-controlled filter system under non-forced, autonomous conditions:
\begin{equation}
\label{MoogVCF}
\dot{\bf x} = {\bf f}\left({\bf x}\right)\qquad\qquad{\rm where}
\qquad\qquad
{\bf x} = \begin{bmatrix} x_{1}\\
x_{2}\\
x_{3}\\
x_{4}\\
\end{bmatrix}\qquad\qquad {\bf f}\left({\bf x}\right) = \omega_{0}\begin{bmatrix} -\tanh(x_{1})-\tanh(\alpha^4 x_{4})\\
-\tanh(x_{2})+\tanh(x_{1})\\
-\tanh(x_{3})+\tanh(x_{2})\\
-\tanh(x_{4})+\tanh(x_{3})\\
\end{bmatrix}
\end{equation}
for constant $\omega_{0}>0$, which serves as a cutoff frequency, and $\alpha = \sqrt{2}r^{1/4}$, for constant $r$ with $0\leq r\leq 1$, which controls the resonance of the filter. Here, ${\bf x}$ is the state, representing the voltages (nondimensionalized) across a series of four capacitors. ${\bf x}$ is a function of time $t$, and $\dot{\bf x}$ indicates the ordinary time derivative of ${\bf x}$. 

What we would like to find for this system is a Lyapunov function $V({\bf x})$, satisfying
\begin{itemize}
\item (C1)\qquad $V({\bf x}) = 0\qquad$ {\rm for} $\qquad {\bf x} = {\bf 0}$
\item (C2)\qquad $V({\bf x}) > 0\qquad$ {\rm for} $\qquad {\bf x} \neq {\bf 0}$
\item (C3)\qquad $\dot{V}\leq 0$ \qquad{\rm for} $\qquad {\bf x} \neq {\bf 0}$ 
\end{itemize}
The third condition above can be strengthened to asymptotic stability using a strict inequality:
\begin{itemize}
\item (C3)'\qquad $\dot{V}< 0$ \qquad{\rm for} $\qquad {\bf x} \neq {\bf 0}$ \qquad{\rm (asymptotic stability)}
\end{itemize}
 If, furthermore, we have the condition of radial unboundedness, or
\begin{itemize}
\item (C4)\qquad $V({\bf x}) \rightarrow \infty\qquad$ {\rm as} $\qquad |{\bf x}| \rightarrow \infty$\qquad{\rm (global asymptotic stability)}
\end{itemize}
we then have a condition for global asymptotic stability. 

We also require that the Lyapunov function satisfy these properties for the whole range of parameters $\omega_{0}>0$ and $0\leq r\leq 1$, though in this short report, we do not assume time-varying behaviour. 

\section{Linearised System}

Under linearisation, the Moog VCF system reduces to
\begin{equation}
\dot{\bf x} = {\bf Ax}\qquad\qquad{\rm where}
\qquad\qquad
{\bf A} = \omega_{0}\begin{bmatrix} -1 & 0 & 0 & -\alpha^4\\
1 & -1 & 0 & 0\\
0 & 1 & -1 & 0\\
0 & 0 & 1 & -1\\
\end{bmatrix}
\end{equation}

For the linearised system, stability is entirely determined by the eigenvalues $\{\lambda_{{\bf A}}\}$ of ${\bf A}$, which are
\begin{equation}
\{\lambda_{{\bf A}}\} = -\omega_{0}+\omega_{0}\alpha\{e^{j\pi/4},\,e^{3j\pi/4},\,e^{5j\pi/4},\,e^{7j\pi/4}\}
\end{equation}
Under the condition that $0\leq r\leq 1$, we then have $0\leq \alpha\leq \sqrt{2}$, and it is true that 
\begin{equation}
{\rm max}\left({\rm Re}\{\lambda_{{\bf A}}\}\right) \leq 0
\end{equation}
and thus all solutions are exponentially non-increasing for the full range of parameters $0\leq r\leq 1$. We would like to determine a Lyapunov function which shows the stability of this system under this same range of parameters. (NB: In this case of a linear system, we could simply diagonalise the system above, but we are interested in extensions to the nonlinear case.) 

The natural choice of a Lyapunov function here is
\begin{equation}
V({\bf x}) = \frac{1}{2}{\bf x}^{T}{\bf x}
\end{equation}
$V$ corresponds to the physical stored energy across the four capacitors in the Moog VCF, and obviously satisfies conditions C1, C2 and C4. It remains to check conditions C3 and C3'. We have:
\begin{equation}
\dot{V} = \dot{{\bf x}}^{T}{\bf x} = {\bf x}^{T}\underbrace{\left(\frac{1}{2}\left({\bf A}+{\bf A}^{T}\right)\right)}_{{\bf A}_{s}}{\bf x}
\end{equation}
However, the matrix ${\bf A}_{s}$ is negative definite only for $0\leq r < 5/12$, for which it satisfies condition C3' (and the weaker condition C3 if $r=5/12$ is included). This choice of Lyapunov function does not recover the full stable parameter range for the linearized Moog VCF. 

A better choice can be achieved through scaling of the state variables. Consider the new variable ${\bf w}$, defined as
\begin{equation}
{\bf w} = {\bf D}{\bf x}\qquad{\rm where}\qquad {\bf D} = \begin{bmatrix} 1 & 0 & 0 & 0\\
0 & \alpha & 0 & 0\\
0 & 0 & \alpha^2 & 0\\
0 & 0 & 0 & \alpha^3\\
\end{bmatrix}
\end{equation}
(Here we leave aside the special case of $\alpha=0$, and assume $0<\alpha\leq \sqrt{2}$.) Now we have
\begin{equation}
\dot{{\bf w}} = {\bf B}{\bf w}\qquad{\rm where}\qquad {\bf B} = \omega_{0}\begin{bmatrix} -1 & 0 & 0 & -\alpha\\
\alpha & -1 & 0 & 0\\
0 & \alpha & -1 & 0\\
0 & 0 & \alpha & -1\\
\end{bmatrix}
\end{equation}
Now choose
\begin{equation}
V({\bf w}) = \frac{1}{2}{\bf w}^{T}{\bf w} = \frac{1}{2}{\bf x}^{T}{\bf D}^2{\bf x}
\end{equation}
which again satisfies C1, C2 and C4. For conditions C3 and C3',
\begin{equation}
\dot{V} = \dot{{\bf w}}^{T}{\bf w} = {\bf w}^{T}\underbrace{\left(\frac{1}{2}\left({\bf B}+{\bf B}^{T}\right)\right)}_{{\bf B}_{s}}{\bf w}
\end{equation}
Now, ${\bf B}_{s}$ is negative definite over the full range of parameters $0<r< 1$ and thus satisfies C3' over this range (and the weaker condition C3 if $r=1$ is included). 

\section{Nonlinear Moog VCF}

Consider now the Moog VCF system as defined in \eqref{MoogVCF}. Define now a general diagonal scaling of the form
\begin{equation}
{\bf w} = {\bf D}{\bf x}\qquad{\rm where}\qquad {\bf D} = \begin{bmatrix} 1 & 0 & 0 & 0\\
0 & d & 0 & 0\\
0 & 0 & d^2 & 0\\
0 & 0 & 0 & d^3\\
\end{bmatrix}
\end{equation}
for some real parameter $d>0$. (Earlier, in the linear case, we set $d=\alpha$, but here we leave $d$ unspecified for the moment.) In terms of ${\bf w}$, the system may now be written as
\begin{equation}
\label{zdef}
\dot{{\bf w}} = \omega_{0}\begin{bmatrix} -\tanh(w_{1})-\tanh(\alpha^4 w_{4}/d^3)\\
-d\tanh(w_{2}/d)+d\tanh(w_{1})\\
-d^2\tanh(w_{3}/d^2)+d^2\tanh(w_{2}/d)\\
-d^3\tanh(w_{4}/d^3)+d^3\tanh(w_{3}/d^2)\\
\end{bmatrix}
=\omega_{0}\underbrace{\begin{bmatrix} -1 & 0 & 0 & -d \\
d & -1 & 0 & 0\\
0 & d & -1 & 0\\
0 & 0 & d & -g\\
\end{bmatrix}}_{{\bf Q}}
\underbrace{\begin{bmatrix} \tanh(w_{1})\\
d\tanh(w_{2}/d)\\
d^2\tanh(w_{3}/d^2)\\
\frac{1}{d}\tanh(\alpha^4 w_{4}/d^3)\\
\end{bmatrix}}_{{\bf z}}
\end{equation}
where note ${\bf z}$ and ${\bf Q}$ as defined above, and where
\begin{equation}
\label{gdef}
g(w_{4}) = d^4 \frac{\tanh(w_{4}/d^3)}{\tanh(\alpha^4 w_{4}/d^3)}
\end{equation}
Notice the $g$ takes on values between $d^4$ and $d^{4}/\alpha^4$. 

\subsection{Lyapunov Function}

Now, define a candidate Lyapunov function for the nonlinear Moog VCF system as
\begin{equation}
\label{Lyapunov_def}
V({\bf w}) = \ln(\cosh(w_{1}))+d^2\ln(\cosh(w_{2}/d))+d^4\ln(\cosh(w_{3}/d^2))+\frac{d^2}{\alpha^4}\ln(\cosh(\alpha^4 w_{4}/d^3))
\end{equation}
(Recall that we have specified $d>0$ and also have restricted $\alpha$ such that $0<\alpha\leq \sqrt{2}$.) $V$ again satisfies conditions C1, C2 and C4 for any choice of $d>0$. For conditions C3 and C3', we have
\begin{equation}
\dot{V} = (\nabla V)^{T}\dot{{\bf w}} = {\bf z}^{T}\dot{{\bf w}}
\end{equation}
using the gradient $\nabla V$ of $V$ with respect to ${\bf w}$, and where ${\bf z}$ is as defined in \eqref{zdef}. But, using the definition of ${\bf Q}$ from \eqref{zdef}, we then have
\begin{equation}
\dot{V} = \omega_{0}{\bf z}^{T}{\bf Q}{\bf z} = \omega_{0}{\bf z}^{T}\underbrace{\left(\frac{1}{2}\left({\bf Q}+{\bf Q}^{T}\right)\right)}_{{\bf Q}_{s}}{\bf z}
\end{equation}
It then remains to determine conditions under which ${\bf Q}_{s}$ is negative definite. To this end, rewrite ${\bf Q}_{s}$ as
\begin{equation}
{\bf Q}_{s} = -{\bf I}_{4}+\frac{d}{2}\underbrace{\begin{bmatrix}0 & 1 & 0 & -1\\
1 & 0 & 1 & 0\\
0 & 1 & 0 & 1\\
-1 & 0 & 1 & f\\
\end{bmatrix}}_{{\bf G}}
\end{equation}
in terms of a matrix ${\bf G}$, which is dependent on $f$, defined in terms of $g$ as
\begin{equation}
f = \frac{2}{d}\left(1-g\right)
\end{equation} 
The maximal positive eigenvalue of ${\bf G}$ may be written as
\begin{equation}
\lambda_{{\bf G}, max} = \max\left(\sqrt{2}, \frac{1}{2}\left(f+\sqrt{f^2+8}\right)\right)
\end{equation}
To satisfy condition C3', we want to find a functional form $d(\alpha)$ such that
\begin{equation}
\label{cond}
\lambda_{{\bf G}, max}<2/d\qquad{\rm for\, all}\qquad \alpha\qquad {\rm with}\qquad 0<\alpha\leq \sqrt{2}
\end{equation}
It is useful at this point to consider two cases for $\alpha$:

\subsubsection*{Case I: $\alpha\leq 1$}
In this case, choose $d=1$. Thus $g$ takes on values between $1/\alpha^4$ and $1$, so $g\geq 1$. This implies that $f\leq 0$, and thus $\lambda_{{\bf G}, max}=\sqrt{2}$, and thus \eqref{cond} is satisfied.
\subsubsection*{Case II: $1<\alpha \leq \sqrt{2}$}
In this case, we can choose $d=\alpha$. Thus $g$ takes on values between $1$ and $\alpha^4$, so $g\geq 1$. Again this implies that $f\leq 0$, and thus $\lambda_{{\bf G}, max}=\sqrt{2}$, and furthermore that 
\begin{equation}
\lambda_{{\bf G}, max} = \sqrt{2} \leq 2/\alpha
\end{equation}
and again \eqref{cond} is satisfied. 

Thus, the choice 
\begin{equation}
d = \max (1,\alpha)
\end{equation}
leads to a Lyapunov function \eqref{Lyapunov_def} that recovers the full range of values of the parameter $\alpha$:
\begin{equation}
\alpha\leq \sqrt{2}\qquad \rightarrow\qquad r\leq 1
\end{equation}



\bibliographystyle{unsrt}
\bibliography{master}

\end{document}